\documentclass[12pt,reqno,a4wide]{amsart}

\usepackage{tikz} 

\usetikzlibrary{decorations.pathreplacing}
\usetikzlibrary{arrows,decorations.pathmorphing,snakes}
\usetikzlibrary{fit}
\usepackage{calrsfs}

\allowdisplaybreaks

\title {Counting edges according to edge-type in $t$-ary trees}

\author{Helmut Prodinger}
\address{Helmut Prodinger\\
	Department of Mathematical Sciences\\
	Stellenbosch University\\
	7602 Stellenbosch,	South Africa  \\ and\\
NITheCS (National Institute for Theoretical and Computational Sciences)\\
South Africa}
\email{hproding@sun.ac.za}


\begin{document}
	
	\begin{abstract}
Using the Lagrange inversion formula, $t$-ary trees are enumerated with respect to edge type (left, middle, right for ternary trees).
	\end{abstract}
	
	\maketitle

	\section{Introduction}
	The  preprint \cite{Burstein} triggered my interest in the sequence A120986 in \cite{OEIS}.  It led to a  consideration of  $T(n,k)$, the
	number of ternary trees with $n$ nodes and $k$ middle edges, see \cite{Three-halves}. 
	
	But we can do better and enumerate such trees with $a_1$ left edges, $a_2$ middle edges, and $a_3$ right edges. Naturally, we must have
	$a_1+a_2+a_3=n-1$. And each permutation of $\{\textit{left}, \textit{middle}, \textit{right} \}$ must lead to the same count, even if different permutations
	are applied on various nodes. We explain everything for ternary trees, but later we will point out how things generalize to $t$-ary trees. We use, apart from the variable $x$ for the nodes, $y_1, y_2, y_3$ for the edge counts; altogether we must look at the coefficient of $x^ny_1^{a_1}y_2^{a_2}y_3^{a_3}$ in a suitable generating function. 
	
	We have 
	\begin{equation*}
G=x(1+y_1G)(1+y_2G)(1+y_3G),
	\end{equation*}
	and, more generally
	\begin{equation*}
		G=x(1+y_1G)\dots(1+y_tG).
	\end{equation*}
	This is already suitable for Lagrange inversion:
	\begin{equation*}
		[x^n]G=\frac{1}{n}[G^{n-1}](1+y_1G)^n\dots(1+y_tG)^n.
	\end{equation*}
	Furthermore
	\begin{align*}
		[x^ny_1^{a_1}\dots y_t^{a_t}]G&=\frac{1}{n}[G^{n-1}y_1^{a_1}\dots y_t^{a_t}](1+y_1G)^n\dots(1+y_tG)^n\\
		&=\frac{1}{n}\binom{n}{a_1}\dots\binom{n}{a_t}.
	\end{align*}

 \section{Wave or particle \emph{aka} lattice path or tree}
 
 Sometimes lattice paths work better than trees, sometimes it is the other way around. 
 Strehl's paths \cite{Strehl-deepest} seem to work better than Kemp's trees \cite{Kemp-deepest}, but the trees in \cite{Naiomi-Kendall} work better that
 the paths in \cite{Naiomi-Kendall}. There are many examples, and the pendulum can go both ways.
 
 One can see from the picture that the trees do sit on different levels, so counting the middle edges cannot be directly justified by counting the residue classes
 $\bmod t$ of the down-steps. However, with our previous argument about the permutations and the large symmetries, even such an approach is feasible. 
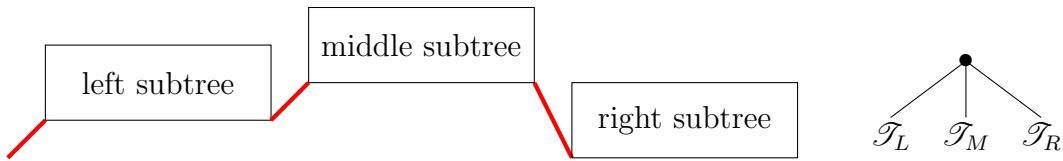
\begin{figure}[h]
	\begin{tikzpicture}[scale=0.5]
\draw[red, ultra thick] (0,0) to (1,1);
\draw (1,1)--(7,1)--(7,3)--(1,3)--(1,1);
	\node at (4,2) {left subtree};
\draw[red, ultra thick] (0+7,1) to (1+7,2);
\draw (1+7,1+1)--(7+7,1+1)--(7+7,3+1)--(1+7,3+1)--(1+7,1+1);
\node at (11,3) {middle subtree};
\draw[red, ultra thick] (14,2) to (15,0);
\draw (1+14,1-1)--(7+14,1-1)--(7+14,3-1)--(1+14,3-1)--(1+14,1-1);
\node at (20-2,1) {right subtree};

	\end{tikzpicture}\qquad
	\begin{tikzpicture}[scale=0.5]
		\node at (0,0) {$\bullet$};
		\node at (-2,-2) {$\mathcal{T}_L$ };
		\draw (0,0)--(-2,-1.5);
		\node at (0,-2) {$\mathcal{T}_M$ };
		\draw (0,0)--(2,-1.5);
		\node at (2,-2) {$\mathcal{T}_R$ };
		\draw (0,0)--(0,-1.5);
	\end{tikzpicture}
	\caption{Ternary tree represented as a lattice path}
	
\end{figure}

\section{Ordered sequences of $t$-ary trees}

We think about $m$ $t$-ary trees, and a giant root (not counted) holding them together. We assume that $1\le m< t$ since then we can count the edges
from the giant root to the regular root of each of the $m$ trees. These edges are labelled by $y_1,\dots,y_m$. Then the Lagrange inversion formula leads to
\begin{equation*}
\frac{m}{n}\binom{n}{a_1-1}\dots \binom{n}{a_m-1}\cdot \binom{n}{a_{m+1}}\dots \binom{n}{a_{t}}.
\end{equation*}

\section{Conclusion}
We hope to report on applications soon.

\end{document}